\documentclass[11pt]{article}
\usepackage{amsmath, graphicx, amsfonts,amssymb, calrsfs}
\usepackage{color}

\def\bbR{\mathbb{R}}
\def\r{\rho}
\def\s{\sigma}
\def\a{\alpha}

\def\dualorlicz{\widetilde{V}_{\phi}}

\def\o{\omega}
\def\ball{B^n_2}

\def\l{\lambda}
\def\cS{\mathcal{S}}

\def\be{\begin{equation}}
\def\ee{\end{equation}}
\def\bea{\begin{eqnarray}}
\def\eea{\end{eqnarray}}
\def\bt{\begin{theorem}}
\def\et{\end{theorem}}
\def\bl{\begin{lemma}}
\def\el{\end{lemma}}
\def\br{\begin{remark}}
\def\er{\end{remark}}
\def\bc{\begin{corollary}}
\def\ec{\end{corollary}}
\def\bd{\begin{definition}}
\def\ed{\end{definition}}
\def\bp{\begin{proposition}}
\def\ep{\end{proposition}}
\newtheorem{theorem}{Theorem}[section]
\newtheorem{lemma}{Lemma}[section]
\newtheorem{remark}{Remark}[section]
\newtheorem{proposition}{Proposition}[section]
\newtheorem{corollary}{Corollary}[section]

\newtheorem{definition}{Definition}[section]
\begin{document}
\title{Dual Orlicz-Brunn-Minkowski theory:   Orlicz $\varphi$-radial addition, Orlicz $L_{\phi}$-dual mixed volume and related inequalities
}

\author{Deping Ye }
\date{ }
\maketitle
 
\begin{abstract}   This paper develops basic setting for the dual Orlicz-Brunn-Minkowski theory for star bodies. An Orlicz $\varphi$-radial addition of two or more star bodies is proposed and related dual Orlicz-Brunn-Minkowski inequality is established. Based on a linear Orlicz $\varphi$-radial addition of two star bodies, we derive a formula for the Orlicz $L_{\phi}$-dual mixed volume. Moreover, a dual Orlicz-Minkowski inequality for the Orlicz $L_{\phi}$-dual mixed volume, a dual Orlicz isoperimetric inequality for the Orlicz $L_{\phi}$-dual surface area and a dual Orlicz-Urysohn inequality for the Orlicz $L_{\phi}$-harmonic mean radius are proved.

 \end{abstract}

\section{Introduction and overview of main results} \label{section introduction}  
  Combination of the Minkowski sum and volume leads to the rich and powerful  classical Brunn-Minkowski theory for convex bodies (i.e., compact convex subsets of $\bbR^n$ with nonempty interiors), which has been thought to be at the core of modern (convex) geometry. Many important results, such as, the Brunn-Minkowski inequality, Minkowski first inequality and the isoperimetric inequality, play fundamental roles in attacking problems in analysis, geometry, quantum information theory, random matrices and many other fields. Readers are referred to the famous book on the classical Brunn-Minkowski theory for convex bodies by Schneider \cite{Sch} for more details and references.  

 In the same spirit, developed by the combination of radial sum and volume, there is a dual theory for star bodies -- the dual  Brunn-Minkowski theory. Such a dual theory was initiated by Lutwak in \cite{Lut1975} and further details were provided in \cite{Lut1988}, where many important objects and fundamental results in the classical Brunn-Minkowski theory have been extended to their dual counterparts. For instance, the dual Minkowski inequality for dual mixed volume is analogous to the Minkowski first inequality for mixed volume, and plays key roles in the solution of the famous Busemann-Petty problem (see e.g., \cite{Gardner1994, GardnerKoldobski1999, Lut1988, Zhang1999}). The literature in the dual Brunn-Minkowski theory is big and continues to grow, see \cite{Bernig2014, GardnerDP2014, Gardner2007, GardnerV1998, GardnerV1999, GardnerV2000, Lutwak1990, Milman2006, Zhang1994} among others. The book written by Gardner \cite{Gardner2005} is an excellent reference for the dual Brunn-Minkowski theory.

 One way to extend the classical Brunn-Minkowski theory and its dual is to replace the linear function $\phi(t)=t$ (note that both Minkowski sum and radial sum are linear) by $\phi(t)=t^p$. With function $\phi(t)=t^p$, the sum for convex bodies is Firey $p$-sum \cite{Firey1962} for $p\geq 1$ (see also \cite{Firey1961} for related work)  and the sum for star bodies is  the $p$-radial sum for $p\neq 0$. Combination of these additions with volume leads to the $L_p$ Brunn-Minkowski theory for convex bodies and its dual (see \cite{Lut1993, Lu1}  and related papers). Readers can find more details and references in Chapter 9 of  \cite{Sch} for the $L_p$ Brunn-Minkowski theory and its dual. 
 
 A further (and more recent) extension of the classical Brunn-Minkowski theory is the new Orlicz-Brunn-Minkowski theory, with homogeneous function $\phi(t)=t^p$ replaced by a nonhomogeneous function $\phi(t)$.  The Orlicz-Brunn-Minkowski theory for convex bodies was initiated from the affine isoperimetric inequalities for Orlicz centroid and projection bodies by Lutwak, Yang and Zhang \cite{LYZ2010a, LYZ2010b}. Note that it is totally nontrivial to find appropriate ways to define Orlicz addition of convex bodies, partially due to lack of homogeneity for nonhomogeneous function $\phi(t)$. Fortunately, this difficulty has been successfully overcome by Gardner, Hug and Weil in their ground breaking paper \cite{Gardner2014}. They gave a nice definition for the Orlicz addition of convex bodies and provided a general framework for the Orlicz-Brunn-Minkowski theory. Based on a linear Orlicz addition for convex bodies, they also derived formulas for the Orlicz $\phi$-mixed volume of two convex bodies. Moreover, they established many important inequalities such as Orlicz-Brunn-Minkowski inequality and Orlicz-Minkowski inequality, whose classical counterparts (namely the  Brunn-Minkowski inequality and Minkowski first inequality) have numerous applications in many fields. Contributions in the  Orlicz-Brunn-Minkowski theory include  \cite{bor2012, bor2013-1, bor2013-2, Chen2011, HaberlLYZ, Li2011, Ludwig2009, Ye2012, Ye2013, Ye2014, Zhu2012} among others.

This paper aims to provide basic setting for the dual Orlicz-Brunn-Minkowski theory for star bodies. In Section \ref{subsection dual brunn minkowski}, the Orlicz $\varphi$-radial addition for two or more star bodies will be introduced and basic properties will be provided. In particular, we will prove a  dual Orlicz-Brunn-Minkowski inequality in  Theorem  \ref{dual brunn minkowski inequality} (see undefined notations in late sections): {\em If $\varphi\in \tilde{\Phi}_2$ and $F_{\varphi}(x_1, x_2)=\varphi(x_1^{-1/n},  x_2^{-1/n})$ is convex, then for $K, L\in \cS_0$,
$$\varphi\bigg(\bigg(\frac{|K\widetilde{+}_{\varphi} L|}{|K|}\bigg)^{1/n}, \bigg(\frac{|K\widetilde{+}_{\varphi} L|}{|L|}\bigg)^{1/n}\bigg)\leq 1.$$ If $F_{\varphi}$ is strictly convex, equality holds if and only if $K$ and $L$ are dilates of each other. }

Section \ref{subsection dual minkowski} is dedicated to develop the formula for the  Orlicz $L_{\phi}$-dual mixed volume of star bodies $K, L$ based on a linear Orlicz $\varphi$-radial addition. Properties of the linear Orlicz $\varphi$-radial addition are similar to those for Orlicz $\varphi$-radial addition proved in Section \ref{subsection dual brunn minkowski}. Related dual Orlicz-Brunn-Minkowski inequality is established in Corollary \ref{dual brunn minkowski inequality--1}. In Section \ref{inequalities}, important inequalities in the classical Brunn-Minkowski theory, such as the Minkowski first inequality, isoperimetric inequality and Urysohn inequality, will be extended to their dual Orlicz counterparts.  For example, in Theorem  \ref{Minkowski inequality -1}, we prove the following dual Orlicz-Minkowski inequality: {\em if $F(t)=\phi(t^{-1/n})$ is convex, then for all $K, L\in \cS_0$, $$\dualorlicz(K, L)\geq |K|\cdot \phi\big({|K|^{1/n}}\cdot {|L|^{-1/n}}\big).$$ If $F(t)$ is strictly convex, equality holds if and only if $K$ and $L$ are dilates of each other.} This dual Orlicz-Minkowski inequality is fundamental in establishing Orlicz affine isoperimetric inequalities for dual Orlicz $L_{\phi}$ affine and geominimal surface areas \cite{Ye2014b}.

 
  \section{Orlicz $\varphi$-radial addition and the dual Orlicz-Brunn-Minkowski inequality}\label{subsection dual brunn minkowski}  
 
 Motivated by the recent elegant work \cite{Gardner2014},  we propose a definition for the Orlicz $\varphi$-radial addition of two or more star bodies in this section. We work on the space $(\bbR^n, \|\cdot\|)$ with $\|\cdot\|$ the usual Euclidean metric induced by the inner product $\langle\cdot, \cdot\rangle$. 
For a measurable set $K\subset \bbR^n$, $|K|$ denotes the Hausdorff content of the appropriate dimension of $K$, in particular, $|K|$ is for the volume of $K$  if $K$ has nonempty interior.  The unit Euclidean ball  in $\bbR^n$ is denoted by $\ball$ and its volume is written by $\o_n$. The spherical measure of the unit sphere  $S^{n-1}$  is denoted by $\s$.  For $\lambda>0$, one defines $\lambda K$ to be  $\lambda K=\{y: y=\lambda x, \ \ x\in K\}. $ The set $K\subset \bbR^n$ is said to be origin-symmetric if $K=-K$ where $-K=\{-x: x\in K\}$. For a linear map $T$, the set $T(K)$ will be written as $TK$ for simplicity.  

In the dual Brunn-Minkowski theory, natural objects are star bodies. A set $K\subset \bbR^n$ is said to be a {\it star body} about the origin, if the line segment from the origin to any point $x\in K$ is contained in $K$ and $K$ has continuous and positive radial function $\r_K(\cdot)$. Here, the {\it radial function}
of $K$, $\r_K: S^{n-1}\rightarrow [0, \infty)$,  is defined by $\r _K(u)=\max \{\l: \l u\in K\}.$  The set of all star bodies (about the origin) in $\bbR^n$ is denoted by  $\cS_0$. Note that $K\in \cS_0$ can be uniquely determined by its radial function $\r_K(\cdot)$ and vice verse.  Two star bodies $K, L\in\cS_0$ are said to be dilates of each other if there is a constant $\lambda>0$ such that $L=\lambda K$, and equivalently $\r_L(u)=\lambda \r_K(u)$ for all $u\in S^{n-1}$. Clearly, for $K, L\in \cS_0$, $$K\subset L\Longleftrightarrow \r_K(u)\leq \r_L(u), \ \ \ \forall u\in S^{n-1}. $$  Moreover, for $K, L\in \cS_0$, one has $$\r_{K\cap L}(u)=\min\{\r_K(u), \r_L(u)\}, \ \ \ \r_{K\cup L}(u)=\max\{\r_K(u), \r_L(u)\}, \ \ \forall u\in S^{n-1}.$$ The natural metric on $\cS_0$ is the radial metric $d_{\r}(\cdot, \cdot): \cS_0\times\cS_0\rightarrow \bbR$ defined as $$d_{\r}(K, L)=\|\r_K-\r_L\|_{\infty}=\sup_{u\in S^{n-1}}|\r_K(u)-\r_L(u)|, \ \ \ \forall K, L\in \cS_0.$$ A sequence of star bodies $\{K_j\}_{j\geq 1}\subset \cS_0$ is said to be convergent to $K\in \cS_0$ in $d_{\r}(\cdot, \cdot)$ if $\lim_{j\rightarrow \infty} d_{\r}(K_j, K)=0$, and equivalently, $\r_{K_j}$ is uniformly convergent to $\r_K$ on $S^{n-1}$.

 Let $m\geq 1$ be a finite integer. Define $\tilde{\Phi}_m$ to be the set of functions $\varphi: [0, \infty)^m \rightarrow [0, \infty)$ such that $\varphi\in \tilde{\Phi}_m$ is continuous, strictly increasing on each variable with $ {\varphi}(0)=0$ and  ${\varphi}(e_i)=1$, where $\{e_i\}_{i=1}^m$ is the standard orthonormal basis in $\bbR^m$.  Let  $\tilde{\Psi}_m$ be the set of functions $\varphi: (0, \infty)^m\rightarrow (0, \infty)$ such that $$\tilde{\varphi}(x_1, \cdots, x_m)=\varphi(1/x_1, \cdots, 1/x_m)\in \tilde{\Phi}_m.$$  Clearly, $\varphi\in \tilde{\Psi}_m$ is continuous and strictly decreasing on each variable. Note that the value of ${\varphi}(e_i)$ is normalized mainly for convenience.  General choices for  ${\varphi}(e_i)\neq 0$ can be taken and our results will follow with only small adjustments.

 For function $\varphi\in \tilde{\Phi}_m$ or  $\varphi\in \tilde{\Psi}_m$, define the Orlicz $\varphi$-radial addition $\widetilde{+}_{\varphi}(K_1, \cdots, K_m)$  of $K_1, \cdots, K_m\in \cS_0$  by the function $\r_{\widetilde{+}_{\varphi} (K_1, \cdots, K_m)}: S^{n-1}\rightarrow \bbR$, such that,  \begin{equation} \varphi \bigg(\frac{\r_{\widetilde{+}_{\varphi}(K_1, \cdots, K_m)}(u)}{\r_{K_1}(u)}, \cdots, \ \frac{\r_{\widetilde{+}_{\varphi}(K_1, \cdots, K_m)}(u)}{\r_{K_m}(u)}\bigg)=1, \ \ \ \forall u\in S^{n-1}. \label{radial harmonic addition}  \end{equation}  For $\varphi(x_1, \cdots, x_m)=x_1^{-p}+\cdots +x_m^{-p}$, one gets the usual $p$-radial addition of $K_1, \cdots, K_m$ (see e.g., \cite{Gardner2002}): $$\big(\r_{\widetilde{+}_{\varphi}(K_1, \cdots, K_m)}(u)\big)^p= \big(\r_{K_1}(u)\big)^p+\cdots +\big(\r_{K_m}(u)\big)^p, \ \ \forall u\in S^{n-1}.$$
When $\varphi(x_1, \cdots, x_m)=x_1^{p}+\cdots +x_m^{p}$ with $p\geq 1$, this is related to the $p$-harmonic combination of $K_1, \cdots, K_m$ (see \cite{Lu1}): $$\big(\r_{\widetilde{+}_{\varphi}(K_1, \cdots, K_m)}(u)\big)^{-p}= \big(\r_{K_1}(u)\big)^{-p}+\cdots +\big(\r_{K_m}(u)\big)^{-p}, \ \ \forall u\in S^{n-1}.$$  In view of these, it seems to be more appropriate to call   $\widetilde{+}_{\varphi}(K_1, \cdots, K_m)$ as the Orlicz $\varphi$-harmonic radial addition (or combination) and the right Orlicz $\varphi$-radial addition should satisfy  \begin{equation*} {\varphi} \bigg(\frac{\r_{K_1}(u)}{\r_{\hat{+}_{\varphi}(K_1, \cdots, K_m)}(u)}, \cdots, \ \frac{\r_{K_m}(u)}{\r_{\hat{+}_{\varphi}(K_1, \cdots, K_m)}(u)}\bigg)=1, \ \ \ \forall u\in S^{n-1}. \end{equation*}  However,  it is easily checked that  $$\widetilde{+}_{\varphi}(K_1, \cdots, K_m)=\hat{+}_{\tilde{\varphi}}(K_1, \cdots, K_m), \ \ with \ \ \tilde{\varphi}(x_1, \cdots, x_m)=\varphi(1/x_1, \cdots, 1/x_m).$$ As equation (\ref{radial harmonic addition}) is relatively easier to analyze and also to be consistent with the function classes $\Phi$ and $\Psi$ defined in \cite{Ye2014} (see also Section \ref{inequalities} in this paper), we will study $\widetilde{+}_{\varphi}(K_1, \cdots, K_m)$ in this paper and call it as the Orlicz $\varphi$-radial addition of $K_1, \cdots, K_m$, instead of the Orlicz $\varphi$-harmonic radial addition. Also note that our definition may be extended to more general functions $\varphi$ and all star shape sets in $\bbR^n$.  For our purpose, we focus on $m=2$ and $K, L\in \cS_0$.     
  
  For  $\varphi\in \tilde{\Phi}_2 $ or  $\varphi \in \tilde{\Psi}_2$, the radial function $\r_{K\widetilde{+}_{\varphi} L}(\cdot)$ for the Orlicz $\varphi$-radial addition of $K, L\in \cS_0$ is determined by the following equation: \begin{equation} \varphi \bigg(\frac{\r_{K\widetilde{+}_{\varphi} L}(u)}{\r_{K}(u)}, \ \frac{\r_{K\widetilde{+}_{\varphi} L}(u)}{\r_{L}(u)}\bigg)=1, \ \ \ \forall u\in S^{n-1}.\label{Orlicz:addition-1} \end{equation} The $p$-radial addition of $K, L\in \cS_0$ is related to function $\varphi(x, y)=x^{-p}+y^{-p}$ and will be denoted by $K\widetilde{+}_{-p} L$, that is, $K\widetilde{+}_{-p} L$ has its radial function $\r_{K\widetilde{+}_{-p} L}(\cdot)$ determined by the equation:   $$ \big(\r_{K}(u)\big)^{p}+\big(\r_{L}(u)\big)^{p}=\big( \r_{K\widetilde{+}_{-p} L}(u)\big)^{p}, \ \ \ \forall u\in S^{n-1}.$$

The following lemma is useful in establishing properties for the Orlicz $\varphi$-radial addition of star bodies. Let $0<a, a_j, b, b_j<\infty$ for all $j\geq 1$ be positive and finite. Let $c, c_j, \tau$ be such that  \begin{equation*} \varphi \bigg(\frac{c_j}{a_j}, \ \frac{c_j}{b_j}\bigg)=1  , \ \   \varphi \left(\frac{c}{a}, \ \frac{c}{b}\right)=1,\  \varphi(\tau, \tau)=1.  \end{equation*}

\bl \label{Property of solution} Let $\varphi\in \tilde{\Phi}_2$ or $\varphi\in \tilde{\Psi}_2$ and real numbers $a, b, c, a_j, b_j, c_j, \tau$ be as above.  \\ (i) $a_1\leq a_2$ and $b_1\leq  b_2$ imply $c_1\leq c_2$. 

\vskip 2mm  \noindent  (ii) Equation $\varphi(x, x)=1$ has $\tau$ as its unique solution. Moreover, $$0<\tau<1 \ \   if \ \ \varphi\in \tilde{\Phi}_2, \ \ and\ \ 1<\tau<\infty\ \ if\ \  \varphi\in \tilde{\Psi}_2.$$  
  \noindent  (iii) Equation  $ \varphi \left( x/{a}, \  x/{b}\right)=1$ also has a unique solution $x=c$ with $$  \tau \cdot \min\{a, b\}\leq c\leq  \tau\cdot  \max\{a, b\}.$$  (iv)  Suppose that $a_j\rightarrow a$ and  $b_j \rightarrow b$ as $j\rightarrow \infty$. Then, $c_j\rightarrow c$ as $j\rightarrow \infty$.   \el

  \noindent{\bf Proof.} (i). Let $a_1\leq  a_2$ and $b_1\leq  b_2$.  Note that $\varphi\in \tilde{\Phi}_2$ is strictly increasing on each variable. Then  $$ \varphi \left(\frac{c_2}{a_2}, \ \frac{c_2}{b_2}\right)=1=\varphi \left(\frac{c_1}{a_1}, \ \frac{c_1}{b_1}\right)\geq  \varphi \left(\frac{c_1}{a_2}, \ \frac{c_1}{b_2}\right).$$ Again by the strictly increasing property of $\varphi$ on each variable, one gets $c_2\geq c_1$.  Similarly, as $\varphi\in \tilde{\Psi}_2$ is strictly decreasing on each variable, one has  $$ \varphi \left(\frac{c_2}{a_2}, \ \frac{c_2}{b_2}\right)=1=\varphi \left(\frac{c_1}{a_1}, \ \frac{c_1}{b_1}\right)\leq \varphi \left(\frac{c_1}{a_2}, \ \frac{c_1}{b_2}\right).$$ Again by the strictly decreasing property of $\varphi$ on each variable, one gets $c_2\geq c_1$. 
  
  \vskip 2mm  \noindent (ii).  Let $\varphi\in \tilde{\Phi}_2$ and $g(x)=\varphi(x, x)$ for all $x\in [0, \infty)$. Then $g(x): [0, \infty)\rightarrow [0, \infty) $ is a continuous, strictly increasing function such that  $g(0)=0$ and $g(1)=\varphi(1,1) > \varphi(0,1)=1$.   Hence, there is one and only one solution for equation $g(x)=1$. This solution must be $\tau$ as $g(\tau)=1$. Clearly $0<\tau<1$ as $g(0)=0$ and $g(1)>1$.   On the other hand, $\tilde{\varphi}(x, y)=\varphi(1/x, 1/y) \in \tilde{\Phi}_2$ for $\varphi\in \tilde{\Psi}_2$. Hence, equation $\varphi(x, x)=\tilde{\varphi}(1/x, 1/x)=1$ has one and only one solution  $\tau$ as $\varphi(\tau, \tau)=1$. Moreover,  $\tilde{\varphi}(1/\tau, 1/\tau)=1$, which leads to $1/\tau \in (0, 1)$ and equivalently $\tau\in (1, \infty)$.
  
   \vskip 2mm  \noindent  (iii).  Let $\varphi\in \tilde{\Phi}_2$. As in the proof of Part (ii), one has $g(x)=\varphi(x/a, x/b)$ is a continuous, strictly increasing function such that  $g(0)=0$ and \begin{eqnarray*} g(\max\{a, b\})&=&\varphi\bigg(\frac{\max\{a, b\}}{a}, \frac{\max\{a, b\}}{b} \bigg)\geq \varphi(1, 1)>\varphi(1, 0)=1. \end{eqnarray*}  Hence, $c$ is the only solution for $g(x)=\varphi(x/a, x/b)=1$ as $g(c)=1$. In particular,  Part (ii) implies that $x=\tau\cdot a$ is the only solution for $\varphi(x/a, x/a)=1$.    Let $a_1=\min\{a, b\}=b_1$. By Part (i), we get $c_1=\tau\cdot \min\{a, b\}\leq c$ as $a_1\leq a$ and $b_1\leq b$. Similarly, by $a\leq \max\{a, b\}$ and $b\leq \max\{a, b\}$, one gets $c\leq \tau\cdot \max\{a, b\}$.  
   
   \vskip 2mm  \noindent   (iv). First of all, if the sequence $\{c_j\}_{j\geq 1}$ has a finite limit, say $c_0$, then $c_0=c$ as $$\varphi(c/a, c/b)=1= \lim_{j\rightarrow \infty} \varphi(c_j/a_j, c_j/b_j)=\varphi(c_0/a, c_0/b).$$ 
  Assume that the sequence $\{c_j\}_{j\geq 1}$ has no limit. Due to $a_j\rightarrow a$ and  $b_j \rightarrow b$ as $j\rightarrow \infty$, there are constants $0<a_0, a_M, b_0, b_M<\infty$ such that $$a_0\leq a_j\leq a_M \ \ and\ \ b_0\leq b_j\leq b_M, \ \ \forall j\geq 1.$$ Part (iii) implies that $$0< \tau \cdot \min\{a_0, b_0\}\leq \tau \cdot \min\{a_j, b_j\} \leq c_j\leq  \tau \cdot \max\{a_j, b_j\}\leq \tau \cdot \max\{a_M, b_M\}<\infty, \ \ \forall j\geq 1. $$ Hence, the sequence $\{c_j\}_{j\geq 1}$ is bounded, and one can always find a convergent subsequence, say $\{c_{j_k}\}_{k\geq 1}$, with limit $\bar{c}\neq c$. Note that $a_{j_k}\rightarrow a$ and $b_{j_k}\rightarrow b$ as $k\rightarrow \infty$. This implies that $\bar{c}$ must be equal to $c$, a contradiction. In conclusion, the sequence $\{c_j\}_{j\geq 1}$ must be convergent to $c$. 
     
\bp \label{Properties: Orlicz radial sum-1} Let $\varphi\in \tilde{\Phi}_2$ or $\varphi\in\tilde{\Psi}_2$, and  $\tau$ be such that $\varphi(\tau, \tau)=1$. \\ (i) If $K, L\in \cS_0$, equation (\ref{Orlicz:addition-1}) uniquely defines a star body  ${K\widetilde{+}_{\varphi} L}\in \cS_0.$\vskip 2mm  \noindent (ii)   For $K, L\in \cS_0$, one has $ \tau \cdot (K\cap L)\subset{K\widetilde{+}_{\varphi} L}\subset \tau\cdot (K\cup L).$ Moreover, ${K\widetilde{+}_{\varphi} L}\subset K\cap L$  if $\phi \in \tilde{\Phi}_2$, while ${K\widetilde{+}_{\varphi} L}\supset K\cup L$ if $\phi \in \tilde{\Psi}_2$. \vskip 2mm  \noindent  (iii) Let $K\in \cS_0$ and $L=\lambda K$ for some $\lambda>0$. Then, ${K\widetilde{+}_{\varphi} L}=\tau_1 \cdot K$ with $\tau_1$ s.t. $\varphi(\tau_1,\tau_1/\lambda)=1$. \vskip 2mm  \noindent  (iv)  Let $K_1\subset K_2$ and $L_1\subset L_2$ be star bodies in $\cS_0$. For all $\varphi\in \tilde{\Phi}_2$ or $\varphi\in \tilde{\Psi}_2$, one has $ K_1\widetilde{+}_{\varphi} L_1 \subset K_2\widetilde{+}_{\varphi} L_2.$
Moreover, if $\varphi\leq \varphi_1$, then \begin{eqnarray*} &&K\widetilde{+}_{\varphi} L\supset K\widetilde{+}_{\varphi_1} L,  \ \ \varphi, \varphi_1\in \tilde{\Phi}_2;\\  &&K\widetilde{+}_{\varphi} L\subset K\widetilde{+}_{\varphi_1} L,  \ \ \varphi, \varphi_1\in \tilde{\Psi}_2.
\end{eqnarray*}   \noindent (v) Let $T\in GL(n)$, the set of all invertible linear transforms from $\bbR^n$ to $\bbR^n$.  For all $K, L\in \cS_0$, one has  $(TK)\widetilde{+}_{\varphi} (TL)=T(K\widetilde{+}_{\varphi} L).$    In particular, $(-K)\widetilde{+}_{\varphi} (-L)=-(K\widetilde{+}_{\varphi} L).$ Moreover, if $K, L\in \cS_0$ are origin-symmetric, then $K\widetilde{+}_{\varphi} L$ is also origin-symmetric. 

\vskip 2mm \noindent (vi) Let $\{K_j\}_{j\geq 1}\subset \cS_0$ and $\{L_j\}_{j\geq 1}\subset \cS_0$ be convergent to $K\in \cS_0$ and $L\in \cS_0$ respectively in $d_{\r}(\cdot, \cdot)$. Then, $ K_j\widetilde{+}_{\varphi} L_j\rightarrow K\widetilde{+}_{\varphi} L$ in $d_{\r}(\cdot, \cdot)$ as $j\rightarrow \infty$.  

\ep

 \noindent {\bf Proof.} (i). Let $u\in S^{n-1}$. Employing Part (iii) of Lemma \ref{Property of solution} to $a=\r_K(u)$ and $b=\r_L(u)$, one gets that $\r_{K\widetilde{+}_{\varphi} L}(u)$ is the unique solution for equation  (\ref{Orlicz:addition-1}) and $$\r_{K\widetilde{+}_{\varphi} L}(u)\geq \tau\cdot \min\{\r_K(u), \r_L(u)\}>0.$$ We claim that  $\r_{K\widetilde{+}_{\varphi} L}(\cdot)$ is continuous at $u\in S^{n-1}$. In fact, for any $u\in S^{n-1}$ given and any sequence $\{u_j\}_{j\geq 1}\subset S^{n-1}$ convergent to $u$, one has for $K, L\in \cS_0$, $$a_j=\r_K(u_j)\rightarrow \r_K(u)=a;\ \ \ b_j=\r_L(u_j)\rightarrow \r_L(u)=b, \ \ \ as\ \ j\rightarrow \infty. $$  Part (iv) of Lemma \ref{Property of solution} implies that $c_j\rightarrow c$ as $j\rightarrow \infty$ with $c_j=\r_{{K\widetilde{+}_{\varphi} L}}(u_j)$ and $c=\r_{K\widetilde{+}_{\varphi} L}(u)$. Thus, $\r_{K\widetilde{+}_{\varphi} L}(\cdot)$ is continuous at $u\in S^{n-1}$ and is continuous on $S^{n-1}$ as desired. In conclusion, equation  (\ref{Orlicz:addition-1}) uniquely determines a star body $K\widetilde{+}_{\varphi} L\in \cS_0$. 
 
 \vskip 2mm \noindent (ii).  Employing Part (iii) of Lemma \ref{Property of solution} to $a=\r_K(u)$ and $b=\r_L(u)$, one gets, for all $u\in S^{n-1}$,  $$\tau \cdot \min\{\r_K(u), \r_L(u)\}\leq \r_{K\widetilde{+}_{\varphi} L}(u) \leq \tau \cdot \max\{\r_K(u), \r_L(u)\}.$$ Equivalently,  $\tau \cdot (K\cap L) \subset  {K\widetilde{+}_{\varphi} L}  \subset \tau \cdot (K\cup L).$  On the other hand, assume that $\r_{K\widetilde{+}_{\varphi} L}(u_0)> \min\{ \r_{K}(u_0), \r_{L}(u_0)\}$ for some $u_0\in S^{n-1}$. Note that $\varphi\in \tilde{\Phi}_2$ is strictly increasing on each variable. Thus, $$1=\varphi\bigg(\frac{\r_{K\widetilde{+}_{\varphi} L}(u_0)}{\r_{K}(u_0)}, \frac{\r_{K\widetilde{+}_{\varphi} L}(u_0)}{\r_{L}(u_0)}\bigg)> \max\{\varphi(0, 1), \varphi(1, 0)\}=1,$$ a contradiction. Hence, for $\varphi\in \tilde{\Phi}_2$, one has $$\r_{K\widetilde{+}_{\varphi} L}(u)\leq \min\{ \r_{K}(u), \r_{L}(u)\}, \ \ \forall u\in S^{n-1} \Longleftrightarrow K\widetilde{+}_{\varphi} L\subset K\cap L.$$ Similarly, for $\varphi\in \tilde{\Psi}_2 $, one has $$\r_{K\widetilde{+}_{\varphi} L}(u)\geq \max\{ \r_{K}(u), \r_{L}(u)\}, \ \ \forall u\in S^{n-1}\Longleftrightarrow K\widetilde{+}_{\varphi} L\supset K\cup L.$$ 
 
 \noindent (iii).  Let $L=\lambda K$ for some $\lambda>0$, which implies $\r_{L}(u)=\lambda \r_K(u)$ for all $u\in S^{n-1}$. One sees that for all $u\in S^{n-1}$,  $\frac{\r_{K\widetilde{+}_{\varphi} L}(u)}{\r_K(u)}$ satisfies the equation $\varphi(x, x/\lambda)=1$ whose unique solution is $\tau_1$. Hence,  $$\frac{\r_{K\widetilde{+}_{\varphi} L}(u)}{\r_K(u)}=\tau_1, \ \  \forall u\in S^{n-1}\ \  \Longleftrightarrow K\widetilde{+}_{\varphi} L=\tau_1\cdot  K.$$
 
\noindent (iv). For $u\in S^{n-1}$, let $a_i=\r_{K_i}(u)$ and $b_i=\r_{L_i}(u)$, $i=1, 2$. Then $a_1\leq a_2$ and $b_1\leq b_2$ due to $K_1\subset K_2$ and $L_1\subset L_2$. Part (i) of Lemma \ref{Property of solution} implies $$c_1=\r_{K_1\widetilde{+}_{\varphi} L_1}(u)\leq c_2=\r_{K_2\widetilde{+}_{\varphi} L_2}(u).$$ Hence, $K_1\widetilde{+}_{\varphi} L_1\subset K_2 \widetilde{+}_{\varphi} L_2$ as desired. 

Now let $\varphi\leq \varphi_1$ with $\varphi, \varphi_1\in \tilde{\Phi}_2$. Then, for all $u\in S^{n-1}$, \begin{eqnarray} \varphi_1 \bigg(\frac{\r_{K\widetilde{+}_{\varphi_1} L}(u)}{\r_{K}(u)}, \ \frac{\r_{K\widetilde{+}_{\varphi_1} L}(u)}{\r_{L}(u)}\bigg)\nonumber &=&1=\varphi \bigg(\frac{\r_{K\widetilde{+}_{\varphi} L}(u)}{\r_{K}(u)}, \ \frac{\r_{K\widetilde{+}_{\varphi} L}(u)}{\r_{L}(u)}\bigg)\\ & \leq& \varphi_1  \bigg(\frac{\r_{K\widetilde{+}_{\varphi} L}(u)}{\r_{K}(u)}, \ \frac{\r_{K\widetilde{+}_{\varphi} L}(u)}{\r_{L}(u)}\bigg).\label{compare varphi-1}\end{eqnarray} As $\varphi_1$ is strictly increasing on each variable, one has, for all $u\in S^{n-1}$,  $$\frac{\r_{K\widetilde{+}_{\varphi_1} L}(u)}{\r_{K}(u)}\leq \frac{\r_{K\widetilde{+}_{\varphi} L}(u)}{\r_{K}(u)}\Longleftrightarrow \r_{K\widetilde{+}_{\varphi_1} L}(u)\leq \r_{K\widetilde{+}_{\varphi} L}(u)\Longleftrightarrow  K\widetilde{+}_{\varphi_1} L \subset  K\widetilde{+}_{\varphi} L. $$ Now let $\varphi\leq \varphi_1$ with $\varphi, \varphi_1\in \tilde{\Psi}_2$. Recall that $\varphi_1$ is strictly decreasing on each variable. By equation (\ref{compare varphi-1}),  one has, for all $u\in S^{n-1}$,  $$\frac{\r_{K\widetilde{+}_{\varphi_1} L}(u)}{\r_{K}(u)}\geq \frac{\r_{K\widetilde{+}_{\varphi} L}(u)}{\r_{K}(u)} \Longleftrightarrow  K\widetilde{+}_{\varphi_1} L \supset  K\widetilde{+}_{\varphi} L. $$  

\noindent (v).  We use $T^{-1}$ for the inverse of $T\in GL(n)$.  For all $u\in S^{n-1}$, let $v=\frac{T^{-1} u}{\|T^{-1} u\|}$ and one has \begin{eqnarray*} \r_{TK}(u)=\sup\{\lambda>0: \lambda u\in TK\}=\sup\{\lambda>0: \lambda T^{-1} u\in K\}=\sup\{\lambda>0: \lambda \|T^{-1} u\| v\in K\}. \end{eqnarray*}  Hence, \begin{equation} \r_{TK}(u)\|T^{-1}u\|=\r_K(v),\label{radial of TK} \end{equation} which further implies  
\begin{eqnarray*}1&=&\varphi \bigg(\frac{\r_{TK\widetilde{+}_{\varphi} TL}(u)}{\r_{TK}(u)}, \ \frac{\r_{TK\widetilde{+}_{\varphi}T L}(u)}{\r_{TL}(u)}\bigg)\\ &=& \varphi \bigg(\frac{\|T^{-1}u\| \r_{TK\widetilde{+}_{\varphi} TL}(u)}{\r_{K}(v)}, \ \frac{\|T^{-1}u\| \r_{TK\widetilde{+}_{\varphi}T L}(u)}{\r_{L}(v)}\bigg)\\ &=&\varphi \bigg(\frac{\r_{K\widetilde{+}_{\varphi} L}(v)}{\r_{K}(v)}, \ \frac{\r_{K\widetilde{+}_{\varphi} L}(v)}{\r_{L}(v)}\bigg). \end{eqnarray*} Part (iii) of Lemma \ref{Property of solution} implies that $\|T^{-1}u\| \r_{TK\widetilde{+}_{\varphi} TL}(u)=\r_{K\widetilde{+}_{\varphi} L}(v)$ and by equation (\ref{radial of TK}) $$(TK)\widetilde{+}_{\varphi} (TL)=T(K\widetilde{+}_{\varphi} L).$$ In particular, $(-K)\widetilde{+}_{\varphi} (-L)=-(K\widetilde{+}_{\varphi} L)$ by letting $T$ be the negative of identity operator. Moreover, if $K=-K$ and $L=-L$, one gets $$-(K\widetilde{+}_{\varphi} L)=(-K)\widetilde{+}_{\varphi} (-L)=K\widetilde{+}_{\varphi} L.$$ Thus, $K\widetilde{+}_{\varphi} L$ is origin-symmetric. 

\vskip 2mm \noindent (vi). Let $\{K_j\}_{j\geq 1}\subset \cS_0$ and $\{L_j\}_{j\geq 1}\subset \cS_0$ be convergent to $K\in \cS_0$ and $L\in \cS_0$ respectively in $d_{\r}(\cdot, \cdot)$. Define $a_j=\r_{K_j}(u)$ and $b_j=\r_{L_j}(u)$ for all $u\in S^{n-1}$. Hence, $a_j\rightarrow a=\r_K(u)$ and $b_j\rightarrow b=\r_L(u)$ as $j\rightarrow \infty$.  By Part (iv) of Lemma \ref{Property of solution}, one gets that for all $u\in S^{n-1}$, $$c_j=\r_{K_j\widetilde{+}_{\varphi} L_j}(u)\rightarrow c=\r_{K\widetilde{+}_{\varphi} L}(u).$$ We now claim that the above convergence is also uniform. To this end, assume that $\r_{K_j\widetilde{+}_{\varphi} L_j}\rightarrow \r_{K\widetilde{+}_{\varphi} L}$ is not uniform on $S^{n-1}$, that is, there is $\varepsilon_0>0$, for all $j\geq 1$, one can find $n_j\geq j$ and $u_{n_j}\in S^{n-1}$, such that, \begin{equation} \label{uniform radial-1--1-1--1--1} |\r_{K_{n_j}\widetilde{+}_{\varphi} L_{n_j}}(u_{n_j})- \r_{K\widetilde{+}_{\varphi} L}(u_{n_j})|\geq \varepsilon_0. \end{equation}  Without loss of generality, assume that $\lim_{j\rightarrow \infty} u_{n_j}=u_0$ for some $u_0\in S^{n-1}$ (as $S^{n-1}$ is compact) and $\lim_{j\rightarrow \infty}\r_{K_{n_j}\widetilde{+}_{\varphi} L_{n_j}}(u_{n_j}) = c_0$  with $c_0$ a positive and finite number. In fact,  as (positive and continuous) $\r_{K_j}(u)$ converges to (positive and continuous) $\r_{K}(u)$ and  (positive and continuous) $\r_{L_j}(u)$ converges to (positive and continuous) $\r_{L}(u)$  uniformly on $S^{n-1}$, one can find $0<M_1<M_2<\infty$, such that, \begin{eqnarray*} M_1&\leq& \r_K(u)\leq M_2, \ \ M_1\leq \r_{K_j}(u)\leq M_2, \ \ \forall u\in S^{n-1}, \ \forall j\geq 1; \\ M_1&\leq& \r_L(u)\leq M_2, \ \ M_1\leq \r_{L_j}(u)\leq M_2, \ \ \forall u\in S^{n-1}, \ \forall j\geq 1.  \end{eqnarray*} Parts (i) and (iii) of Lemma \ref{Property of solution}  imply that \begin{eqnarray*} \tau M_1 \leq \r_{K\widetilde{+}_{\varphi} L}(u)\leq \tau M_2, \ \ \tau M_1\leq \r_{K_{j}\widetilde{+}_{\varphi} L_{j}}(u)\leq \tau M_2, \ \ \forall u\in S^{n-1}, \ \forall j\geq 1.  \end{eqnarray*} Hence,  the sequence  $\r_{K_{n_j}\widetilde{+}_{\varphi} L_{n_j}}(u_{n_j})$ is bounded and has a convergent subsequence. 

In conclusion, we have $\lim_{j\rightarrow \infty} u_{n_j}=u_0$  and $\lim_{j\rightarrow \infty} \r_{K_{n_j}\widetilde{+}_{\varphi} L_{n_j}}(u_{n_j})=c_0$. This further implies \begin{eqnarray*}1=\varphi \bigg(\frac{\r_{K_{n_j}\widetilde{+}_{\varphi} L_{n_j}}(u_{n_j})}{\r_{K_{n_j}}(u_{n_j})}, \ \frac{\r_{K_{n_j}\widetilde{+}_{\varphi} L_{n_j}}(u_{n_j})}{\r_{L_{n_j}}(u_{n_j})}\bigg)\rightarrow  \varphi \bigg(\frac{c_0}{\r_{K}(u_0)}, \ \frac{c_0}{\r_{L}(u_0)}\bigg) \end{eqnarray*}  where we have used the uniform convergence of $\r_{K_j}(u)\rightarrow \r_K(u)$ and $\r_{L_j}(u)\rightarrow \r_L(u)$.  This implies that $c_0=\r_{K\widetilde{+}_{\varphi} L}(u_0)$. On the other hand, equation (\ref{uniform radial-1--1-1--1--1}) implies \begin{equation*} \label{uniform radial-1--1-1--1} |\r_{K_{n_j}\widetilde{+}_{\varphi} L_{n_j}}(u_{n_j})- \r_{K\widetilde{+}_{\varphi} L}(u_{n_j})|\rightarrow |c_0-\r_{K\widetilde{+}_{\varphi} L}(u_0)|\geq \varepsilon_0, \ \ \ as \ j\rightarrow \infty,\end{equation*} a contradiction with  $c_0=\r_{K\widetilde{+}_{\varphi} L}(u_0)$. Hence,  $\r_{K_j\widetilde{+}_{\varphi} L_j}\rightarrow \r_{K\widetilde{+}_{\varphi} L}$ uniformly on $S^{n-1}$, and equivalently, $K_j\widetilde{+}_{\varphi} L_j\rightarrow K\widetilde{+}_{\varphi} L$ in $d_{\r}(\cdot, \cdot)$ as $j\rightarrow \infty$.

\vskip 2mm The following theorem provides a dual Orlicz-Brunn-Minkowski inequality for $|K\widetilde{+}_{\varphi} L|.$
\bt \label{dual brunn minkowski inequality} {\bf (Dual Orlicz-Brunn-Minkowski inequality).} Let $\varphi\in \tilde{\Phi}_2$ or $\varphi\in \tilde{\Psi}_2$, and  $F_{\varphi}(x_1, x_2)=\varphi(x_1^{-1/n},  x_2^{-1/n})$.  \vskip 2mm \noindent (i) If $F_{\varphi}$ is convex, then for $K, L\in \cS_0$,
$$\varphi\bigg(\bigg(\frac{|K\widetilde{+}_{\varphi} L|}{|K|}\bigg)^{1/n}, \bigg(\frac{|K\widetilde{+}_{\varphi} L|}{|L|}\bigg)^{1/n}\bigg)\leq 1.$$ If in addition $F_{\varphi}$ is strictly convex, equality holds if and only if $K$ and $L$ are dilates of each other. \vskip 2mm \noindent (ii) If $F_{\varphi}$ is concave, then for $K, L\in \cS_0$,
$$\varphi\bigg(\bigg(\frac{|K\widetilde{+}_{\varphi} L|}{|K|}\bigg)^{1/n}, \bigg(\frac{|K\widetilde{+}_{\varphi} L|}{|L|}\bigg)^{1/n}\bigg)\geq 1.$$ If in addition $F_{\varphi}$ is strictly concave,  equality holds if and only if $K$ and $L$ are dilates of each other.   \et

\noindent{\bf Proof.}
(i). Jensen's inequality  and equation (\ref{Orlicz:addition-1}) imply that \begin{eqnarray*} 1 &=& \int_{S^{n-1}}\varphi \bigg(\frac{\r_{K\widetilde{+}_{\varphi} L}(u)}{\r_{K}(u)}, \ \frac{\r_{K\widetilde{+}_{\varphi} L}(u)}{\r_{L}(u)}\bigg) \frac{[\r_{K\widetilde{+}_{\varphi} L}(u)]^n }{n |K\widetilde{+}_{\varphi} L|}\,d\s(u) \\ &=& \int_{S^{n-1}}F_{\varphi} \bigg(\frac{[\r_{K}(u)]^n}{[\r_{K\widetilde{+}_{\varphi} L}(u)]^n}, \ \frac{[\r_{L}(u)]^n}{[\r_{K\widetilde{+}_{\varphi} L}(u)]^n}\bigg)\ \frac{[\r_{K\widetilde{+}_{\varphi} L}(u)]^n }{n |K\widetilde{+}_{\varphi} L|}\,d\s(u) \\ &\geq& F_{\varphi} \bigg( \int_{S^{n-1}}\frac{[\r_{K}(u)]^n }{n |K\widetilde{+}_{\varphi} L|} \,d\s(u), \int_{S^{n-1}}\frac{[\r_{L}(u)]^n }{n |K\widetilde{+}_{\varphi} L|} \,d\s(u)\bigg) \\ &=&F_{\varphi} \bigg(\frac{|K| }{ |K\widetilde{+}_{\varphi} L|}, \ \ \frac{|L|}{ |K\widetilde{+}_{\varphi} L|}  \bigg) \\ &=& \varphi\bigg(\bigg(\frac{|K\widetilde{+}_{\varphi} L|}{|K|}\bigg)^{1/n}, \bigg(\frac{|K\widetilde{+}_{\varphi} L|}{|L|}\bigg)^{1/n}\bigg).\end{eqnarray*}  If in addition $F_{\varphi}$ is strictly convex, then equality holds if and only if   there are two constants $c_1, c_2>0$ (note that $K, L, K\widetilde{+}_{\varphi} L\in \cS_0$ have continuous and positive radial functions), such that, $$\frac{[\r_{K}(u)]^n}{[\r_{K\widetilde{+}_{\varphi} L}(u)]^n}=c_1^n \ \& \  \frac{[\r_{L}(u)]^n}{[\r_{K\widetilde{+}_{\varphi} L}(u)]^n}=c_2^n\ \ \Longrightarrow \r_L(u)=\frac{c_2}{c_1}\r_K(u), \ \ \  \forall   u\in S^{n-1}.$$ That is,  $K$ and $L$ are dilates of each other.  
\vskip 2mm \noindent (ii). Let $F_{\varphi}$ be concave, then \begin{eqnarray*} 1&=&   \int_{S^{n-1}}F_{\varphi} \bigg(\frac{[\r_{K}(u)]^n}{[\r_{K\widetilde{+}_{\varphi} L}(u)]^n}, \ \frac{[\r_{L}(u)]^n}{[\r_{K\widetilde{+}_{\varphi} L}(u)]^n}\bigg) \frac{[\r_{K\widetilde{+}_{\varphi} L}(u)]^n }{n |K\widetilde{+}_{\varphi} L|}\,d\s(u) \\ &\leq& F_{\varphi} \bigg( \int_{S^{n-1}}\frac{[\r_{K}(u)]^n }{n |K\widetilde{+}_{\varphi} L|} \,d\s(u), \int_{S^{n-1}}\frac{[\r_{L}(u)]^n }{n |K\widetilde{+}_{\varphi} L|} \,d\s(u)\bigg) \\ &=& \varphi\bigg(\bigg(\frac{|K\widetilde{+}_{\varphi} L|}{|K|}\bigg)^{1/n}, \bigg(\frac{|K\widetilde{+}_{\varphi} L|}{|L|}\bigg)^{1/n}\bigg).\end{eqnarray*}  Similar to case (i),  if in addition $F_{\varphi}$ is strictly concave, equality holds in the dual Orlicz-Brunn-Minkowski inequality if and only if $K, L\in \cS_0$ are dilates of each other. \\

\noindent{\bf Remark.} For $\varphi(x, y)=x^{-p}+y^{-p}$, our dual Orlicz-Brunn-Minkowski inequality becomes the well-known $L_p$-dual Brunn-Minkowski inequality: Let $K, L\in \cS_0$, then for $p\in (n, \infty)\cup (-\infty, 0)$, $$|K\widetilde{+}_{-p} L|^{p/n}\geq |K|^{p/n}+|L|^{p/n};$$ while for $p\in  (0, n]$, $$|K\widetilde{+}_{-p} L|^{p/n}\leq |K|^{p/n}+|L|^{p/n}.$$ In particular, the case $p=1$ is called the dual Brunn-Minkowski inequality and the case $p=n-1$ is called the dual Kneser-S\"{u}ss inequality. See the survey on the Brunn-Minkowski inequality by Gardner \cite{Gardner2002}.

\section{The linear Orlicz $\varphi$-radial addition and the Orlicz $L_{\phi}$-dual  mixed volume} \label{subsection dual minkowski}

For our purpose to derive formulas for the Orlicz $L_{\phi}$-dual  mixed volume, we need the {\it linear Orlicz $\varphi$-radial addition}.  Let $\varphi(x, y)=\a\phi_1(x)+\beta\phi_2(y)$ and $\varphi_i(x, y)=\a_i\phi_1(x)+\beta_i\phi_2(y)$ where functions $\phi_1, \phi_2$  are either both in $\tilde{\Phi}_1$ or both in $\tilde{\Psi}_1$, and $\a, \beta, \a_i, \beta_i>0$ are constants such that $\a+\beta\geq 1$ and $\a_i+\beta_i\geq 1$ for $i=1, 2$. Again,  the assumption $\a+\beta\geq 1$ is mainly for convenience and general choices can be taken.

\bd\label{Linear Orlicz---1-1} For $K, L\in \cS_0$,  let $K\widetilde{+}_{\varphi}L$ be determined by the radial function $\r_{K\widetilde{+}_{\varphi} L}(u)$ s.t. \begin{equation}\label{Linear orlicz--1} 1=\a \phi_1\bigg(\frac{\r_{K\widetilde{+}_{\varphi} L}(u)}{\r_{K}(u)} \bigg)+\beta\phi_2\bigg(\frac{\r_{K\widetilde{+}_{\varphi} L}(u)}{\r_{L}(u)}\bigg), \ \ \ \forall u\in S^{n-1}. \end{equation}  \ed
\noindent {\bf Remark.} Similar to the proof of Proposition \ref{Orlicz:addition-1}, we can see that $K\widetilde{+}_{\varphi} L$ is a star body uniquely determined by equation (\ref{Linear orlicz--1}).  If $K$ and $L$ are dilates of each other, say $L=\lambda K$ for some $\lambda>0$, then $K\widetilde{+}_{\varphi} L=\tau_1 K$ with $\tau_1$ such that \begin{equation}\label{Dilates of each other--1} \a \phi_1(\tau_1)+\beta \phi_2(\tau_1/\lambda)=1. \end{equation} Moreover, as $\a+\beta\geq 1$, there is a unique and finite number $\tau>0$ s.t. $\a \phi_1(\tau)+\beta\phi_2(\tau)=1$. Hence,  $$\tau (K\cap L)\subset K\widetilde{+}_{\varphi} L\subset \tau (K\cup L). $$ If $\a_1\leq \a_2$ and $\beta_1\leq \beta_2$, then for $\phi_1, \phi_2\in \tilde{\Phi}_1$, \begin{eqnarray*} 1&=&  \a_2 \phi_1\bigg(\frac{\r_{K\widetilde{+}_{ \varphi_2} L}(u)}{\r_{K}(u)} \bigg)+\beta_2 \phi_2\bigg(\frac{\r_{K\widetilde{+}_{\varphi_2} L}(u)}{\r_{L}(u)}\bigg) \\ &=&\a_1 \phi_1\bigg(\frac{\r_{K\widetilde{+}_{\varphi_1} L}(u)}{\r_{K}(u)} \bigg)+\beta_1 \phi_2\bigg(\frac{\r_{K\widetilde{+}_{\varphi_1} L}(u)}{\r_{L}(u)}\bigg)\\ &\leq&  \a_2 \phi_1\bigg(\frac{\r_{K\widetilde{+}_{\varphi_1} L}(u)}{\r_{K}(u)} \bigg)+\beta_2 \phi_2\bigg(\frac{\r_{K\widetilde{+}_{\varphi_1} L}(u)}{\r_{L}(u)}\bigg), \ \ \ \ \forall u\in S^{n-1}.\end{eqnarray*}  Due to the strictly increasing property of $\phi_1, \phi_2\in \tilde{\Phi}_1$, one has $$\r_{K\widetilde{+}_{\varphi_2} L}(u)\leq \r_{K\widetilde{+}_{\varphi_1} L}(u), \ \ \ \forall u\in S^{n-1} \Longleftrightarrow K\widetilde{+}_{\varphi_2} L \subset  K\widetilde{+}_{\varphi_1} L.$$  Similarly, if $\a_1\leq \a_2$ and $\beta_1\leq \beta_2$, then for $\phi_1, \phi_2\in \tilde{\Psi}_1$ (and hence $\phi_1, \phi_2$ are strictly decreasing), $$\r_{K\widetilde{+}_{\varphi_2} L}(u)\geq \r_{K\widetilde{+}_{\varphi_1} L}(u), \ \ \ \forall u\in S^{n-1} \Longleftrightarrow  K\widetilde{+}_{\varphi_2} L \supset  K\widetilde{+}_{\varphi_1} L.$$ 

The following corollary is the related dual Orlicz-Brunn-Minkowski inequality for the linear Orlicz $\varphi$-radial addition.  For functions $\phi_1, \phi_2$, let $F_1(t)= \phi_1(t^{-1/n})$ and $F_2(t)=\phi_2(t^{-1/n})$. Recall that $\varphi(x, y)=\a \phi_1(x)+\beta \phi_2(y)$.

\bc \label{dual brunn minkowski inequality--1} {\bf (Dual Orlicz-Brunn-Minkowski inequality).}  Let  $\phi_1, \phi_2$ be either both in $\tilde{\Phi}_1$ or both in  $\tilde{\Psi}_1$.  \vskip 2mm \noindent (i)  Let $K, L\in \cS_0$. If both $F_1(t)$ and $F_2(t)$ are convex,  one has  $$\a \phi_1\bigg(\bigg(\frac{|K\widetilde{+}_{\varphi} L|}{|K|}\bigg)^{1/n}\bigg)+\beta \phi_2\bigg(\bigg(\frac{|K\widetilde{+}_{\varphi} L|}{|L|}\bigg)^{1/n}\bigg)\leq 1.$$ If in addition at least one of $F_1(t)$ and $F_2(t)$  are strictly convex, equality holds if and only if $K$ and  $L$ are dilates of each other. \vskip 2mm \noindent (ii)  Let $K, L\in \cS_0$. If both $F_1(t)$ and $F_2(t)$ are concave,  one has
$$\a\phi_1\bigg(\bigg(\frac{|K\widetilde{+}_{\varphi} L|}{|K|}\bigg)^{1/n}\bigg)+\beta\phi_2\bigg(\bigg(\frac{|K\widetilde{+}_{\varphi} L|}{|L|}\bigg)^{1/n}\bigg)\geq 1.$$ If in addition at least one of $F_1(t)$ and $F_2(t)$ are  strictly concave, equality holds if and only if $K$ and  $L$ are dilates of each other.  \ec
\noindent{\bf Proof.} The proof is similar to that of Theorem \ref{dual brunn minkowski inequality}. For completeness, we include a brief proof here. First, note that equation (\ref{Linear orlicz--1})  is equivalent to $$\a F_{1}\bigg(\frac{[\r_{K}(u)]^n}{[\r_{K\widetilde{+}_{\varphi} L}(u)]^n}\bigg)+\beta F_2\bigg(\frac{[\r_{L}(u)]^n}{[\r_{K\widetilde{+}_{\varphi} L}(u)]^n}\bigg)=1, \ \ \ \forall u\in S^{n-1}.$$  (i).  Let  $F_1(t)$ and $F_2(t)$ be both convex. Jensen's inequality implies that 
 \begin{eqnarray} 1   &=& \int_{S^{n-1}} \bigg[ \a F_{1}\bigg(\frac{[\r_{K}(u)]^n}{[\r_{K\widetilde{+}_{\varphi} L}(u)]^n}\bigg)+\beta F_2\bigg(\frac{[\r_{L}(u)]^n}{[\r_{K\widetilde{+}_{\varphi} L}(u)]^n}\bigg)\bigg]  \frac{[\r_{K\widetilde{+}_{\varphi} L}(u)]^n }{n |K\widetilde{+}_{\varphi} L|}\,d\s(u)\nonumber \\ &\geq& \a F_{1} \bigg( \int_{S^{n-1}}\frac{[\r_{K}(u)]^n }{n |K\widetilde{+}_{\varphi} L|} \,d\s(u)\bigg)+\beta F_2\bigg(\int_{S^{n-1}}\frac{[\r_{L}(u)]^n }{n |K\widetilde{+}_{\varphi} L|} \,d\s(u)\bigg) \nonumber\\   &=& \a \phi_1\bigg(\bigg(\frac{|K\widetilde{+}_{\varphi} L|}{|K|}\bigg)^{1/n}\bigg)+\beta \phi_2\bigg(\bigg(\frac{|K\widetilde{+}_{\varphi} L|}{|L|}\bigg)^{1/n}\bigg).\label{strictly equal--1}\end{eqnarray}  Clearly, by equation (\ref{Dilates of each other--1}), if $K$ and $L$ are dilates of each other, one gets the equality.

  Assume that in addition at least one of $F_1(t)$ and $F_2(t)$  are strictly convex. Without loss of generality, let $F_1(t)$ be strictly convex. Then, equality holds in inequality (\ref{strictly equal--1}) if and only if equalities hold for Jensen's inequality on both $F_1(t)$ and $F_2(t)$. As $F_1(t)$ is strictly convex,  there is a constant $c_1>0$ (note that $K,  K\widetilde{+}_{\varphi} L\in \cS_0$ have continuous positive radial functions), such that,  $$\frac{[\r_{K}(u)]^n}{[\r_{K\widetilde{+}_{\varphi} L}(u)]^n}=\bigg(\frac{1}{c_1}\bigg)^n, \ \ \ \forall u\in S^{n-1} \ \Longleftrightarrow \r_{K\widetilde{+}_{\varphi} L}(u)=c_1\cdot \r_K(u), \ \ \  \forall   u\in S^{n-1}.$$   Together with equation (\ref{Linear orlicz--1}), one has, for all $u\in S^{n-1}$,  \begin{equation*} 1=\a \phi_1(c_1)+\beta\phi_2\bigg(\frac{c_1\cdot \r_K(u)}{\r_{L}(u)}\bigg)\ \Longrightarrow\  \frac{c_1\cdot \r_K(u)}{\r_{L}(u)}=\phi_2^{-1}\bigg(\frac{1-\a \phi_1(c_1)}{\beta}\bigg). \end{equation*} That is, $K$ and $L$ are dilates of each other as desired.

\vskip 2mm \noindent (ii). Let  $F_1(t)$ and $F_2(t)$ be both concave. Jensen's inequality implies that 
 \begin{eqnarray*} 1   &=& \int_{S^{n-1}} \bigg[ \a F_{1} \bigg(\frac{[\r_{K}(u)]^n}{[\r_{K\widetilde{+}_{\varphi} L}(u)]^n}\bigg)+\beta F_2\bigg(\frac{[\r_{L}(u)]^n}{[\r_{K\widetilde{+}_{\varphi} L}(u)]^n}\bigg)\bigg]  \frac{[\r_{K\widetilde{+}_{\varphi} L}(u)]^n }{n |K\widetilde{+}_{\varphi} L|}\,d\s(u) \\ &\leq& \a F_{1} \bigg( \int_{S^{n-1}}\frac{[\r_{K}(u)]^n }{n |K\widetilde{+}_{\varphi} L|} \,d\s(u)\bigg)+\beta F_2\bigg(\int_{S^{n-1}}\frac{[\r_{L}(u)]^n }{n |K\widetilde{+}_{\varphi} L|} \,d\s(u)\bigg) \\   &=& \a \phi_1\bigg(\bigg(\frac{|K\widetilde{+}_{\varphi} L|}{|K|}\bigg)^{1/n}\bigg)+\beta \phi_2\bigg(\bigg(\frac{|K\widetilde{+}_{\varphi} L|}{|L|}\bigg)^{1/n}\bigg).\end{eqnarray*} Similar to Part (i), if in addition at least one of $F_1(t)$ and $F_2(t)$ are  strictly concave, equality holds if and only if $K$ and  $L$ are dilates of each other. \\

The case $\a=1$ and $\beta=\epsilon$ is important for us to derive the formula for Orlicz $L_{\phi}$-dual mixed volume. Let $K\widetilde{+}_{\epsilon,\phi_1, \phi_2} L\in \cS_0$ be determined by the radial function $\r_{K\widetilde{+}_{\epsilon,\phi_1, \phi_2} L}(\cdot)$, such that, \begin{equation*} 1=\phi_1\bigg(\frac{\r_{K\widetilde{+}_{\epsilon, \phi_1, \phi_2} L}(u)}{\r_{K}(u)} \bigg)+\epsilon\phi_2\bigg(\frac{\r_{K\widetilde{+}_{\epsilon, \phi_1, \phi_2} L}(u)}{\r_{L}(u)}\bigg), \ \ \ \forall u\in S^{n-1}. \end{equation*} By the remark after Definition \ref{Linear Orlicz---1-1}, one gets, for all $0<\epsilon\leq 1$,\begin{eqnarray} \label{linear compare-1} &&K\widetilde{+}_{1, \phi_1, \phi_2} L\subset K\widetilde{+}_{\epsilon, \phi_1, \phi_2} L\subset K, \ \ \ \phi_1, \phi_2\in \tilde{\Phi}_1; \\    \label{linear compare-2} &&K\widetilde{+}_{1, \phi_1, \phi_2} L\supset K\widetilde{+}_{\epsilon, \phi_1, \phi_2} L\supset K, \ \ \ \phi_1, \phi_2\in \tilde{\Psi}_1.\end{eqnarray}

 A more convenient equivalent formula is 
\begin{equation} \label{definition:linear combination} 1=G_1\bigg( \frac{[\r_{K\widetilde{+}_{\epsilon, \phi_1, \phi_2} L}(u)]^n}{[\r_{K}(u)]^n}\bigg)+\epsilon G_2\bigg( \frac{[\r_{K\widetilde{+}_{\epsilon, \phi_1, \phi_2} L}(u)]^n}{[\r_{L}(u)]^n}\bigg), \ \ \ \forall u\in S^{n-1}, \end{equation} where $G_1(t)=\phi_1(t^{1/n})$ and $G_2(t)=\phi_2(t^{1/n})$. 
 
\bl\label{uniform convergence} Let $K, L\in \cS_0$. \\ (i) Let $\phi_1, \phi_2$ be either both in $\tilde{\Phi}_1$ or both in $\tilde{\Psi}_1$. The following limit is uniform on $S^{n-1}$ $$\r_{K\widetilde{+}_{\epsilon,\phi_1, \phi_2} L}(u)\rightarrow \r_K(u)\ \ as\ \ \epsilon\rightarrow 0^+.$$ (ii)    Let $\phi_1, \phi_2\in \tilde{\Phi}_1$ be such that $\phi_{1, l}'(1)$, the left-derivative of  $\phi_1(t)$ at $1$, exists and is finite. The following limit is uniform on $S^{n-1}$  $$ \phi_{1, l}'(1)  \lim_{\epsilon\rightarrow 0^+} \frac{[\r_K(u)]^n-[\r_{K\widetilde{+}_{\epsilon,\phi_1, \phi_2} L}(u)]^n}{ n  \epsilon}= \phi_2\bigg(\frac{\r_K(u)}{\r_L(u)}\bigg) [\r_K(u)]^n.$$  (iii) Let $\phi_1, \phi_2\in \tilde{\Psi}_1$ be such that $\phi_{1, r}'(1)$, the right-derivative of  $\phi_1(t)$ at $1$, exists and is finite. The following limit is uniform on $S^{n-1}$   $$\phi_{1, r}'(1)   \lim_{\epsilon\rightarrow 0^+} \frac{[\r_K(u)]^n-[\r_{K\widetilde{+}_{\epsilon,\phi_1, \phi_2} L}(u)]^n}{n \epsilon}= \phi_2\bigg(\frac{\r_K(u)}{\r_L(u)}\bigg) [\r_K(u)]^n.$$  \el 

\noindent {\bf Proof.} (i). Let $\phi_1, \phi_2\in \tilde{\Phi}_1$, formula (\ref{linear compare-1}) and the monotone increasing property of $\phi_1, \phi_2$ (and hence $\phi_1^{-1}$) imply that \begin{eqnarray*} \phi_1^{-1}\left(1-\epsilon\phi_2\left(a\right)\right)   \leq   \frac{\r_{K\widetilde{+}_{\epsilon, \phi_1, \phi_2} L}(u)}{\r_{K}(u)} = \phi_1^{-1}\bigg(1-\epsilon\phi_2\bigg(\frac{\r_{K\widetilde{+}_{\epsilon, \phi_1, \phi_2} L}(u)}{\r_{L}(u)}\bigg)\bigg) \nonumber \leq  \phi_1^{-1}\left(1-\epsilon\phi_2\left(b\right)\right),  \end{eqnarray*} uniformly on $S^{n-1}$ with $0<a, b<\infty$ defined by (note that  $K, L, K\widetilde{+}_{1, \phi_1, \phi_2} L\in \cS_0$ have continuous positive radial functions) $$a=\max_{u\in S^{n-1}} \frac{\r_K(u)}{\r_L(u)}, \ \ b=\min _{u\in S^{n-1}} \frac{\r_{K\widetilde{+}_{1, \phi_1, \phi_2} L}(u)}{\r_L(u)} .$$ Taking $\epsilon\rightarrow 0^+$, we get that the desired limit is uniform on $S^{n-1}$. The case $\phi_1, \phi_2\in \tilde{\Psi}_1$ follows by a similar argument.

\vskip 2mm \noindent (ii). Note that $G_{1, l}'(1)=\phi_{1, l}'(1)/n$. It is enough to prove that $$ G_{1, l}'(1)  \lim_{\epsilon\rightarrow 0^+} \frac{1-\big(\frac{\r_{K\widetilde{+}_{\epsilon,\phi_1, \phi_2} L}(u)}{\r_K(u)}\big)^n}{ \epsilon}= \phi_2\left(\frac{\r_K(u)}{\r_L(u)}\right).$$  As $\phi_1\in \tilde{\Phi}_1$ is a strictly increasing function, $G_1(t)$ is also increasing and  its inverse exists (denoted by $G_1^{-1}(t)$). Given $u\in S^{n-1}$, let  $$z(\epsilon)=  \bigg(\frac{\r_{K\widetilde{+}_{\epsilon,\phi_1, \phi_2} L}(u)}{\r_K(u)}\bigg)^n.$$ Equation (\ref{definition:linear combination}) implies that \begin{eqnarray*} \frac{1- \big(\frac{\r_{K\widetilde{+}_{\epsilon,\phi_1, \phi_2} L}(u)}{\r_K(u)}\big)^n}{\epsilon} &=& \frac{1- z(\epsilon)}{\epsilon}   = \phi_2\bigg( \frac{\r_{K\widetilde{+}_{\epsilon, \phi_1, \phi_2} L}(u)}{\r_{L}(u)}\bigg) \cdot \frac{ 1- z(\epsilon)}{1-G_1(z(\epsilon))}. \end{eqnarray*}  If $\epsilon \rightarrow 0^+$, then $z(\epsilon)\rightarrow 1^-$ uniformly on $S^{n-1}$ by Part (i) and inequality (\ref{linear compare-1}).  Therefore, the following limit is uniformly on $S^{n-1}$, \begin{eqnarray*} \lim_{\epsilon\rightarrow 0^+} \frac{1- \big(\frac{\r_{K\widetilde{+}_{\epsilon,\phi_1, \phi_2} L}(u)}{\r_K(u)}\big)^n}{\epsilon} &=& \lim_{ \epsilon \rightarrow 0^+}\phi_2\bigg( \frac{\r_{K\widetilde{+}_{\epsilon, \phi_1, \phi_2} L}(u)}{\r_{L}(u)}\bigg) \cdot \lim_{z(\epsilon)\rightarrow 1^-} \frac{ 1- z(\epsilon)}{1-G_1(z(\epsilon))}\\    &=& \frac{1}{G_{1, l}'(1)}\cdot  \phi_2\bigg(\frac{\r_K(u)}{\r_L(u)}\bigg)= \frac{n}{\phi_{1, l}'(1)}\cdot  \phi_2\bigg(\frac{\r_K(u)}{\r_L(u)}\bigg). \end{eqnarray*}

\noindent  (iii). Note that $G_{1, r}'(1)=\phi_{1, r}'(1)/n$.  For $\phi_1, \phi_2\in \tilde{\Psi}_1$, Part (i) and inequality (\ref{linear compare-2}) imply that $z(\epsilon)\rightarrow 1^+$ uniformly on $S^{n-1}$ as $\epsilon\rightarrow 0^+$. Therefore, uniformly on $S^{n-1}$, \begin{eqnarray*} \lim_{\epsilon\rightarrow 0^+} \frac{1- \big(\frac{\r_{K\widetilde{+}_{\epsilon,\phi_1, \phi_2} L}(u)}{\r_K(u)}\big)^n}{\epsilon} &=&\lim_{\epsilon \rightarrow 0^+} \phi_2\bigg( \frac{\r_{K\widetilde{+}_{\epsilon, \phi_1, \phi_2} L}(u)}{\r_{L}(u)}\bigg) \cdot \lim_{z(\epsilon)\rightarrow 1^+}  \frac{ 1- z(\epsilon)}{1-G_1(z(\epsilon))}\\    &=& \frac{n}{\phi_{1, r}'(1)}\cdot  \phi_2\bigg(\frac{\r_K(u)}{\r_L(u)}\bigg). \end{eqnarray*}  

The asymptotic behaviour of $|K\widetilde{+}_{\epsilon,\phi_1, \phi_2} L|$ is stated in the following theorem. See similar results for $\phi_1(t)=\phi_2(t)=t^p$ with $p\geq 1$ in \cite{Lu1}.  

\bt \label{limit of dual mixed volume} Let $K, L\in \cS_0$. \vskip 2mm \noindent (i) Let $\phi_1, \phi_2\in \tilde{\Phi}_1$ be such that $\phi_{1, l}'(1)$ exists and is finite.  One has $$ \phi_{1, l}'(1)  \lim_{\epsilon\rightarrow 0^+} \frac{|K|-|K\widetilde{+}_{\epsilon,\phi_1, \phi_2} L|}{n\epsilon}=\frac{1}{n}\int_{S^{n-1}} \phi_2\bigg(\frac{\r_K(u)}{\r_L(u)}\bigg) [\r_K(u)]^n \,d\s(u).$$ (ii) Let $\phi_1, \phi_2\in \tilde{\Psi}_1$ be such that $\phi_{1, r}'(1)$ exists and is finite.  One has   $$ \phi_{1, r}'(1)  \lim_{\epsilon\rightarrow 0^+} \frac{|K|- |K\widetilde{+}_{\epsilon,\phi_1, \phi_2} L|}{n\epsilon}=\frac{1}{n} \int_{S^{n-1}} \phi_2\left(\frac{\r_K(u)}{\r_L(u)}\right) [\r_K(u)]^n \,d\s(u).$$ \et 

\noindent {\bf Proof.} Part (ii) of Lemma \ref{uniform convergence} implies that for $\phi_1, \phi_2\in \tilde{\Phi}_1$ with $\phi_{1, l}'(1)$ finite,  \begin{eqnarray*}  \phi_{1, l}'(1)  \lim_{\epsilon\rightarrow 0^+} \frac{|K|-|K\widetilde{+}_{\epsilon,\phi_1, \phi_2} L|}{n\epsilon}  &=&  \phi_{1, l}'(1)  \lim_{\epsilon\rightarrow 0^+} \int_{S^{n-1}} \frac{[\r_K(u)]^n-[\r_{K\widetilde{+}_{\epsilon,\phi_1, \phi_2} L}(u)]^n}{n^2 \epsilon} \,d\s(u)  \\ &=&\phi_{1, l}'(1)  \int_{S^{n-1}} \bigg(\lim_{\epsilon\rightarrow 0^+} \frac{ [\r_K(u)]^n-[\r_{K\widetilde{+}_{\epsilon,\phi_1, \phi_2} L}(u)]^n   }{n^2 \epsilon}\bigg) \,d\s(u) \\ &=&\frac{1}{n}\int_{S^{n-1}} \phi_2\left(\frac{\r_K(u)}{\r_L(u)}\right) [\r_K(u)]^n \,d\s(u).\end{eqnarray*}   
By Part (iii) of Lemma   \ref{uniform convergence}, one can prove the  case $\phi_1, \phi_2\in \tilde{\Psi}_1$ along the same line.  

In view of Theorem \ref{limit of dual mixed volume}, we define the Orlicz $L_{\phi}$-dual mixed volume of $K, L\in \cS_0$ for all continuous positive functions $\phi$ as follows.  
\bd\label{dual mixed volume:definition} Let $\phi: (0, \infty)\rightarrow (0, \infty)$ be a continuous positive function. Define the Orlicz $L_{\phi}$-dual  mixed volume, denoted by $\dualorlicz(K, L)$,  of star bodies $K, L\in \cS_0$ as  $$\dualorlicz(K, L)=\frac{1}{n}\int _{S^{n-1}}\phi\bigg(\frac{\r_K(u)}{\r_L(u)}\bigg) [\r_K(u)]^n \,d\s(u).$$ \ed When $\phi(t)=1$, one gets $\dualorlicz(K, L)=|K|$ for all $K, L\in \cS_0$. If $\phi(t)=t^{-n}$, one has $\dualorlicz(K, L)=|L|$ for all $K, L\in \cS_0$.

Denote by $SL(n)$ the subset of $GL(n)$ with unit absolute value of determinant, i.e., $T\in SL(n)$ if $T\in GL(n)$ with $|det(T)|=1$. We now prove that the Orlicz $L_{\phi}$-dual  mixed volume $\dualorlicz(K, L)$ is $SL(n)$-invariant. 
 \bt\label{affine invariance dual orlicz mixed volume} For $T\in SL(n)$ and $K, L\in \cS_0$, one has  
 \begin{eqnarray*}
 \dualorlicz(TK, TL)= \dualorlicz(K, L). \end{eqnarray*}   \et
 \noindent {\bf Proof.}   Let $T\in SL(n)$. Define $u=\frac{T^{-1} v}{\|T^{-1} v\|}\in S^{n-1}$ for $v\in S^{n-1}$. Equation (\ref{radial of TK})  implies that $$\phi\left(\frac{\r_K(u)}{\r_L(u)}\right)=\phi\left(\frac{\r_{TK}(v)}{\r_{TL}(v)}\right).$$
  Note that $\frac{1}{n} [\r_K(u)]^n\,d\s(u)$ is the volume element of $K$ and hence  $ [\r_{TK}(v)]^n\,d\s(v)=[\r_K(u)]^n\,d\s(u)$. Thus, for $ T\in SL(n)$ and all $K, L\in \cS_0$, 
\begin{eqnarray*}
n\dualorlicz(TK, TL)&=& \int _{S^{n-1}} \!\phi \bigg(\!\frac{\r_{TK}(v)}{\r_{TL}(v)}\!\bigg) [\r_{TK}(v)]^n\,d\s(v)\\ &=&  \!\! \int _{S^{n-1}}\! \phi\bigg(\!\frac{\r_{K}(u)}{\r_{L}(u)}\!\bigg) [\r_{K}(u)]^n\,d\s(u) \!=\! n  \dualorlicz(K, L). \end{eqnarray*}

 Let $\psi(t)$ have inverse function $\psi^{-1}(t)$ and $H(t)=(\phi \circ \psi^{-1})(t)$ be the composition function of $\phi(t)$ and $\psi^{-1}(t)$. The following result compares $\dualorlicz(K, L)$ and $\widetilde{V}_{\psi}(K, L)$ for all $K, L\in \cS_0$.  
 \bp Let $K, L\in \cS_0$ and $\phi, \psi$ be continuous positive functions with $H(t)=(\phi \circ \psi^{-1})(t)$. \vskip 2mm \noindent 
 (i) If $H(t)$ is convex, then $$\frac{\dualorlicz(K, L)}{|K|} \geq H\bigg(\frac{\widetilde{V}_{\psi}(K, L)}{|K|}\bigg).$$  If in addition $H(t)$ is strictly convex, equality holds if and only if $K$ and $L$ are dilates of each other. \vskip 2mm \noindent (ii) If $H(t)$ is concave, then $$\frac{\dualorlicz(K, L)}{|K|} \leq H\bigg(\frac{\widetilde{V}_{\psi}(K, L)}{|K|}\bigg).$$ If in addition $H(t)$ is strictly concave, equality holds if and only if $K$ and $L$ are dilates of each other.\ep
 
 \noindent {\bf Proof.} (i). Let  $H(t)$ be convex. Jensen's inequality implies that  \begin{eqnarray*}\frac{\dualorlicz(K, L)}{|K|}&=& \frac{1}{n|K|}\int_{S^{n-1}} \phi\bigg(\frac{\r_K(u)}{\r_L(u)}\bigg) [\r_K(u)]^n\,d\s(u)\\ &=& \frac{1}{n|K|}\int_{S^{n-1}} H\bigg(\psi\bigg(\frac{\r_K(u)}{\r_L(u)}\bigg)\bigg) [\r_K(u)]^n\,d\s(u)\\&  \geq&  H\bigg(  \frac{1}{n|K|}\int_{S^{n-1}} \psi\bigg(\frac{\r_K(u)}{\r_L(u)}\bigg) [\r_K(u)]^n\,d\s(u)\bigg) = H\bigg(\frac{\widetilde{V}_{\psi}(K, L)}{|K|}\bigg). \end{eqnarray*} Note that $K, L\in \cS_0$ have continuous positive radial functions.  Therefore, if $H(t)$ is strictly convex, equality holds in the above inequality if and only if there is a constant $\lambda>0$, such that  $$\psi\bigg(\frac{\r_K(u)}{\r_L(u)}\bigg)=\psi(\lambda)\ \ \Longrightarrow\ \  \frac{\r_K(u)}{\r_L(u)}=\lambda, \ \ \ \ \ \ \forall u\in S^{n-1}.$$ That is $K$ and $L$ are dilates of each other. 
 
\vskip 2mm \noindent (ii). Let $H(t)$ be concave. Jensen's inequality implies that  \begin{eqnarray*}\frac{\dualorlicz(K, L)}{|K|} = \frac{1}{n|K|}\int_{S^{n-1}} H\bigg(\psi\bigg(\frac{\r_K(u)}{\r_L(u)}\bigg)\bigg) [\r_K(u)]^n\,d\s(u)\leq H\bigg(\frac{\widetilde{V}_{\psi}(K, L)}{|K|}\bigg). \end{eqnarray*} Similar to Part (i), if in addition $H(t)$ is strictly concave, equality holds if and only if $K$ and $L$ are dilates of each other.

 \section{Inequalities for the Orlicz $L_{\phi}$-dual mixed volume}\label{inequalities}
 
 This section dedicates to generalize fundamental inequalities in the classical Brunn-Minkowski theory for convex bodies, such as, Minkowski first inequality, isoperimetric inequality and Urysohn inequality,  to their dual Orlicz counterparts. 

For  $\phi: (0, \infty)\rightarrow (0, \infty)$, let $F(t)=\phi(t^{-1/n})$ and hence $\phi(t)=F(t^{-n})$. Consider the sets of functions $\Phi$ and $\Psi$ as \begin{eqnarray*} \Phi\!\!\!&=&\!\!\! \{\mbox{$\phi: (0, \infty)\rightarrow (0, \infty)$ :  $F(t)$ is either a constant or a convex function}\}, \\  \Psi\!\!\!&=&\!\!\!\{\phi: (0, \infty)\rightarrow (0, \infty): F(t)\ \mbox{is either a constant or an increasing  concave function}\}.\end{eqnarray*}   Sample functions in $\Phi$ are: $t^p$ with $p\in (-\infty, -n]\cup (0, \infty)$ and convex increasing functions; sample  functions in $\Psi$ are: $t^p$ with $p\in [-n, 0)$,  $\arctan(t^{-n})$ and $\ln (1+t^{-n})$. Both $\Phi$ and $\Psi$ may not contain some nice functions such as $\phi(t)=e^{-t}$ and could have neither convex nor concave functions. Note that functions $\phi(t)$ and $F(t)$ have opposite monotonicity: if one is increasing then the other will be decreasing (and vice verse).  The relation of convexity and concavity between $\phi(t)$ and $F(t)$ is not clear; however if one is convex and increasing, then  the other must be convex and decreasing. See \cite{Ye2014} for more discussion on the properties of $\Phi$ and $\Psi$. (Note that $\Phi$ and $\Psi$ defined  in \cite{Ye2014} have extra restrictions: strict convexity and strict concavity respectively).   
 
\bt \label{Minkowski inequality -1} {\bf (Dual Orlicz-Minkowski inequality).} Let $K, L\in \cS_0$. \\ (i) For $\phi\in \Phi$, one has $$\dualorlicz(K, L)\geq |K|\cdot \phi\big({|K|^{1/n}}\cdot {|L|^{-1/n}}\big).$$ If in addition $F(t)$ is strictly convex, equality holds if and only if $K$ and $L$ are dilates of each other. \vskip 2mm   \noindent (ii) For $\phi\in \Psi$, one has $$\dualorlicz(K, L)\leq |K|\cdot \phi\big({|K|^{1/n}}\cdot {|L|^{-1/n}}\big).$$ If in addition $F(t)$ is strictly concave, equality holds if and only if $K$ and $L$ are dilates of each other.   \et 
\noindent {\bf Proof.} (i). Let $\phi\in \Phi$ and then $F(t)=\phi(t^{-1/n})$ is convex.  Jensen's inequality (for convex function $F(t)$) implies that for all $K, L\in {\cS_0}$,  \begin{eqnarray}  \dualorlicz(K, L) &=&\frac{1}{n} \int_{S^{n-1}} \phi\left(\frac{\r_K(u)}{ \r_L(u)} \right) [\r_K(u)]^n\,d\s(u) \nonumber\\  &=&   |K| \int_{S^{n-1}} F\bigg(\frac{ \r_L^n(u)}{\r_K^n (u)} \bigg) \frac{[\r_K(u)]^n}{n|K|}  \,d\s(u)\nonumber \\  &\geq &  |K| \cdot F\bigg( \int_{S^{n-1}} \frac{ \r_L^n(u)} {n|K|}  \,d\s(u)\bigg) \nonumber \\  &=&  |K|\cdot  F\bigg(\frac{|L|}{ |K|}  \bigg) =  \phi\big({|K|^{1/n}}\cdot {|L|^{-1/n}}) \cdot  |K|. \label{isoperimetric:dual--1----1}
 \end{eqnarray} If in addition $F$ is strictly convex, then equality holds in inequality (\ref{isoperimetric:dual--1----1}) if and only if there is $\lambda>0$, such that $\r_L(u)=\lambda \r_K(u)$ for all $u\in S^{n-1}$.  That is, $K$ and $L$ are dilates of each other.
 
  \vskip 2mm  \noindent 
 (ii). Let $\phi\in \Psi$ and then $F(t)=\phi(t^{-1/n})$ is concave.  Jensen's inequality (for concave function $F(t)$) implies that for all $K, L\in {\cS_0}$,  \begin{eqnarray}  \dualorlicz(K, L)    &=&  |K| \int_{S^{n-1}} F\left(\frac{ \r_L^n(u)}{\r_K^n (u)} \right) \frac{[\r_K(u)]^n}{n|K|}  \,d\s(u)\nonumber \\  &\leq &  |K| \cdot F\bigg( \int_{S^{n-1}} \frac{ \r_L^n(u)} {n|K|}  \,d\s(u)\bigg)  =  \phi\big({|K|^{1/n}}\cdot {|L|^{-1/n}}) \cdot  |K|. \label{isoperimetric:dual--2----2}
 \end{eqnarray} If in addition $F$ is strictly concave, then equality holds in inequality (\ref{isoperimetric:dual--2----2}) if and only if there is $\lambda>0$, such that, $\r_L(u)=\lambda \r_K(u)$ for all $u\in S^{n-1}$. That is, $K$ and $L$ are dilates of each other.

Define the  Orlicz $L_\phi$-dual surface area of $K$ to be $n\dualorlicz (K, \ball)$ and denote by $\widetilde{S}_{\phi}(K)$. It is easily checked that for all $r>0$, $$\widetilde{S}_{\phi}(r \ball)=n\dualorlicz (r\ball, \ball)= \phi(r)\cdot n|r\ball|. $$ We use  $\widetilde{S}_p(K)$ for the case  $\phi(t)=t^p$, and  for $\lambda>0$, \begin{equation} \label{homogeneous p dual} \widetilde{S}_p(\lambda K)=\lambda^{n+p}\widetilde{S}_p(K).\end{equation} 

We now establish the following dual Orlicz isoperimetric inequality for $\widetilde{S}_{\phi}(K)$. Let $B_K$ be the origin-symmetric Euclidean ball with $|B_K|=|K|.$ Therefore, $B_K=r\ball $ with $r={|K|^{1/n}}{|\ball|^{-1/n}},$ and then \begin{equation} \label{dual surface area of BK} \widetilde{S}_{\phi}(B_K)= \phi(r)\cdot n|r\ball|= \phi\bigg(\frac{|K|^{1/n}}{|\ball|^{1/n}}\bigg)\cdot n|K|. \end{equation}  
 
\bt {\bf (Dual Orlicz isoperimetric inequality).} Let $K\in \cS_0$.  
\vskip 2mm \noindent  (i) For $\phi\in \Phi$, one has $$\widetilde{S}_{\phi}(K)\geq \widetilde{S}_{\phi}(B_K).$$ If in addition $F(t)$ is strictly convex, equality holds if and only if $K$ is an origin-symmetric Euclidean ball. \vskip 2mm \noindent  (ii) For $\phi\in \Psi$, one has $$\widetilde{S}_{\phi}(K)\leq \widetilde{S}_{\phi}(B_K).$$ If in addition $F(t)$ is strictly concave, equality holds if and only if $K$ is an origin-symmetric Euclidean ball. \et

 \noindent {\bf Proof.} (i). Let $K\in \cS_0$.    The dual Orlicz-Minkowski inequality implies that for $\phi\in \Phi$, one has   $$\widetilde{S}_{\phi}(K)=n\dualorlicz(K, \ball)\geq n|K|\cdot \phi\big({|K|^{1/n}}\cdot {|\ball|^{-1/n}}\big)=\widetilde{S}_{\phi}(B_K)$$ where the last equality follows from formula (\ref{dual surface area of BK}). If in addition $F(t)$ is strictly convex, equality holds if and only if $K$ and $\ball$ are dilates of each other.  That is, $K$ is an origin-symmetric Euclidean ball.  
 
\vskip 2mm  \noindent (ii). For $\phi\in \Psi$, the dual Orlicz-Minkowski inequality implies that  $$\widetilde{S}_{\phi}(K)=n\dualorlicz(K, \ball)\leq n|K|\cdot \phi\big({|K|^{1/n}}\cdot {|\ball|^{-1/n}}\big)=\widetilde{S}_{\phi}(B_K)$$ where the last equality follows from formula (\ref{dual surface area of BK}).  If in addition $F(t)$ is strictly concave, equality holds if and only if $K$ and $\ball$ are dilates of each other.  That is, $K$ is an origin-symmetric Euclidean ball.  
  
 \vskip 2mm \noindent {\bf Remark.} The dual Orlicz isoperimetric inequality states that for all star bodies with fixed volume, the origin-symmetric Euclidean ball has the minimal Orlicz $L_\phi$-dual surface area for $\phi\in \Phi$. If $\phi(t)=t^p$ with $p\in (0, \infty)\cup (-\infty, -n)$, one can even have, by formula (\ref{homogeneous p dual}), $$\frac{\widetilde{S}_p(K)}{\widetilde{S}_p(\ball)}\geq \bigg(\frac{|K|}{|\ball|}\bigg)^{\frac{n+p}{n}}. $$ On the other hand, if  $\phi\in \Psi$, for all star bodies with fixed volume, the origin-symmetric Euclidean ball has the maximal Orlicz $L_\phi$-dual surface area. If $\phi(t)=t^p$ with $p\in (-n, 0)$, one can even have, by formula (\ref{homogeneous p dual}), $$\frac{\widetilde{S}_p(K)}{\widetilde{S}_p(\ball)}\leq \bigg(\frac{|K|}{|\ball|}\bigg)^{\frac{n+p}{n}}. $$
 
For $K\in \cS_0$, define the Orlicz $L_{\phi}$-harmonic mean radius of $K$ (denoted by $\widetilde{\o}_{\phi}(K)$)  as $$\widetilde{\o}_{\phi}(K)=\frac{1}{n\o_n} \int_{S^{n-1}}\phi\bigg(\frac{1}{\r_K(u)}\bigg)\,d\s(u)=\frac{\dualorlicz(\ball, K)}{\o_n}.$$ Clearly, for $K=r\ball$, then $$\widetilde{\o}_{\phi}(r\ball)=\frac{1}{n\o_n} \int_{S^{n-1}}\phi\big(1/r\big)\,d\s(u)=\phi\big(1/r\big).$$ In particular, as $B_K=r\ball$ with $r=|K|^{1/n}|\ball|^{-1/n}$, one has $$\widetilde{\o}_{\phi}(B_K) =\phi\big({|\ball|^{1/n}}\cdot {|K|^{-1/n}} \big).$$ One can also define the  Orlicz $L_{\phi}$ mean radius of $K$ as  $$\frac{1}{n\o_n} \int_{S^{n-1}}\phi\big({\r_K(u)}\big)\,d\s(u)=\widetilde{\o}_{\tilde{\phi}}(K)=\frac{1}{n\o_n} \int_{S^{n-1}}\tilde{\phi}\bigg(\frac{1}{\r_K(u)}\bigg)\,d\s(u)$$ where $\tilde{\phi}(t)=\phi(t^{-1})$. To be consistent with function classes $\Phi$ and $\Psi$, we prove the following dual Orlicz-Urysohn inequality for the Orlicz $L_{\phi}$-harmonic mean radius of $K$. 
\bt {\bf (Dual Orlicz-Urysohn inequality).} Let $K\in \cS_0$. \vskip 2mm \noindent (i) If $\phi\in \Phi$, then $$\widetilde{\o}_{\phi}(K)\geq \widetilde{\o}_{\phi}(B_K).$$  If in addition $F(t)$ is strictly convex, equality holds if and only if  $K$ is an origin-symmetric Euclidean ball. \vskip 2mm \noindent 
(ii) If $\phi\in \Psi$, then $$\widetilde{\o}_{\phi}(K)\leq \widetilde{\o}_{\phi}(B_K).$$  If in addition $F(t)$ is strictly concave, equality holds if and only if $K$ is an origin-symmetric Euclidean ball. 
 \et \noindent {\bf Proof.} (i). The dual Orlicz-Minkowski inequality implies that for $\phi\in \Phi$, $$\dualorlicz(\ball, K)\geq |\ball|\cdot \phi\big({|\ball|^{1/n}}\cdot {|K|^{-1/n}}\big).$$ Dividing both sides by $\o_n=|\ball|$, one gets, 
 $$\widetilde{\o}_{\phi}(K)=\frac{\dualorlicz(\ball, K)}{\o_n}\geq   \phi\big({|\ball|^{1/n}}\cdot {|K|^{-1/n}}\big)
 =\widetilde{\o}_{\phi}(B_K).$$ If in addition $F(t)$ is strictly convex, equality holds if and only if $K$ and $\ball$ are dilates of each other. That is, $K$ is an origin-symmetric Euclidean ball.  
 
 \vskip 2mm \noindent (ii). The dual Orlicz-Minkowski inequality implies that for $\phi\in \Psi$, $$\dualorlicz(\ball, K)\leq |\ball|\cdot \phi\big({|\ball|^{1/n}}\cdot {|K|^{-1/n}}\big).$$ Dividing both sides by $\o_n=|\ball|$, one gets, 
 $$\widetilde{\o}_{\phi}(K)=\frac{\dualorlicz(\ball, K)}{\o_n}\leq   \phi\big({|\ball|^{1/n}}\cdot {|K|^{-1/n}}\big)=
 \widetilde{\o}_{\phi}(B_K).$$  If in addition $F(t)$ is strictly concave, equality holds if and only if  $K$ and $\ball$ are dilates of each other. That is, $K$ is an origin-symmetric Euclidean ball.

\vskip 2mm \noindent {\bf Acknowledgments.} The research of DY is supported
 by a NSERC grant.

\vskip 2mm \noindent Deping Ye, \ \ \ {\small \tt deping.ye@mun.ca}\\
{\small \em Department of Mathematics and Statistics\\
   Memorial University of Newfoundland\\
   St. John's, Newfoundland, Canada A1C 5S7 }

\end{document}